\documentclass{amsart}
\usepackage{amsmath,amssymb}

\begin{document}


\newcommand{\R}{   \ensuremath{ \mathbb{R} }}
\newcommand{\C}{   \ensuremath{ \mathbb{C} }}
\newcommand{\N}{   \ensuremath{ \mathbb{N} }}
\newcommand{\Z}{   \ensuremath{ \mathbb{Z} }}
\newcommand{\Q}{   \ensuremath{ \mathbb{Q} }}
\newcommand{\F}{   \ensuremath{ \mathbb{F} }}


\newcommand{\CH}{   \ensuremath{ \mathcal{H}  }}   
\newcommand{\V}{    \ensuremath{ \mathcal{V}  }}   
\newcommand{\A}{    \ensuremath{ \mathcal{A}  }}   
\newcommand{\W}{    \ensuremath{ \mathcal{W}  }}   
\newcommand{\D}{    \ensuremath{ \mathcal{D}  }}   
\newcommand{\I}{    \ensuremath{ \mathfrak{I} }}   
\newcommand{\PS}{   \ensuremath{ \mathfrak{P} }}   
\newcommand{\Aut}{  \ensuremath{ \mathrm{Aut} }}   


\newcommand{\zz}{z,\bar{z}}
\newcommand{\zzs}{z,\bar{z},s}

\newcommand{\as}{      \ensuremath{ \mbox{  as  } }}
\newcommand{\unifas}{  \ensuremath{ \mbox{  uniformly as  } }}
\newcommand{\st}{      \ensuremath{ \mbox{  s.~t.  } }}


\renewcommand{\Im}{\ensuremath{ \mathrm{Im} }}
\renewcommand{\Re}{\ensuremath{ \mathrm{Re} }}

\renewcommand{\bar}[1]{   \ensuremath{ \overline{#1} }}
\renewcommand{\hat}[1]{   \widehat{ #1 }}
\renewcommand{\tilde}[1]{ \widetilde{ #1 }}


\newcommand{\dxys}{\, dx \, dy \, ds }
\newcommand{\app}[2]{\ensuremath{ \int_{\R^{#1}} {#2} \dxys }}


\newcommand{\ddth}{       \ensuremath{ \left. \frac{d}{d \theta} \right|_0 }}
\newcommand{\SC}{         \ensuremath{ \mathcal{S} \hat{\mathbf{C}} }}
\newcommand{\SB}{         \ensuremath{ \mathcal{S} \mathbf{b} }}
\newcommand{\vpc}{        \ensuremath{ \nu_A \rho_A \hat{C}|_{S^1} }}
\newcommand{\covec}[1]{   \ensuremath{ {#1}^T \, dZ + \bar{#1}^T \, d \bar{Z} }}


\newcommand{\Arg}{    \ensuremath{ \mathrm{Arg} }}
\newcommand{\bv}{     \ensuremath{ \mathrm{bv}\, }}
\newcommand{\codim}{  \ensuremath{ \mathrm{codim} }}
\newcommand{\coef}{   \ensuremath{ \mathrm{coef} }}
\newcommand{\defect}{ \ensuremath{ \mathrm{defect} }}
\newcommand{\deriv}{  \ensuremath{ \mathrm{deriv} }}
\newcommand{\diam}{   \ensuremath{ \mathrm{diam} }}
\newcommand{\disk}{   \ensuremath{ \mathrm{disk} }}
\newcommand{\eval}{   \ensuremath{ \mathrm{eval} }}
\newcommand{\del}{    \partial }
\newcommand{\Hol}{    \ensuremath{ \mathrm{Hol} }}
\newcommand{\holO}{   \ensuremath{ \mathcal{O}}}
\newcommand{\id}{     \ensuremath{ \mathrm{id}}}
\newcommand{\inc}{    \ensuremath{ \mathrm{inc}}}
\newcommand{\image}{  \ensuremath{ \mathrm{image}\, }}
\newcommand{\Levi}{   \ensuremath{ \mathcal{L} }}
\newcommand{\maxid}{  \ensuremath{ \mathfrak{m} }}
\newcommand{\Null}{   \ensuremath{ \mathrm{Null} }}
\newcommand{\Orb}{    \ensuremath{ \mathrm{Orb} }}
\newcommand{\Pow}{    \ensuremath{ \mathcal{P} }}
\newcommand{\psub}{   \ensuremath{ \stackrel{\textstyle \subset}{\neq} }}
\newcommand{\sgn}{    \ensuremath{ \mathrm{sgn} }}
\newcommand{\Span}{   \ensuremath{ \mathrm{span} }}
\newcommand{\supp}{   \ensuremath{ \mathrm{supp}\, }}
\newcommand{\rank}{   \ensuremath{ \mathrm{rank} }}
\newcommand{\tensor}{ \otimes}


\newcommand{\conv}[3]{      \ensuremath{ {#1} \to {#2}  \mbox{  in  } {#3} }}
\newcommand{\inner}[1]{     \ensuremath{ \langle {#1} \rangle }}
\newcommand{\norm}[1]{      \ensuremath{ \left\| #1 \right\| }}
\newcommand{\setof}[2]{     \ensuremath{ \left\{ #1 : #2 \right\} }}
\newcommand{\setofbar}[2]{  \ensuremath{ \left\{ #1 \, \left| \, #2 \right. \right\} }}


\newcommand{\dd}[2]{    \ensuremath{ \frac{ \displaystyle \del {#1} }{ \displaystyle \del {#2} } }}
\newcommand{\ddat}[2]{  \ensuremath{ \left. \frac{ \displaystyle \del  }{ \displaystyle \del {#1} } \right|_{#2} }}


\newcommand{\Jac}[2]{   \ensuremath{ \left( \dd{#1}{#2} \right) }}
\newcommand{\Jacat}[3]{ \ensuremath{ \left( \dd{#1}{#2} \left( {#3} \right) \right) }}


\newcommand{\Strut}{\rule[-4mm]{0mm}{10mm}}  

\newcommand{\zero}{        \ensuremath{  \mathbf{ 0 }  }}
\newcommand{\zeromatrix}{  \ensuremath{  \mbox{\LARGE\textbf{0}}  }}
\newcommand{\ident}{       \ensuremath{  \mbox{\LARGE\textbf{$I$}} }}


\newtheorem{thm}{Theorem}[section]
\newtheorem{lem}[thm]{Lemma}
\newtheorem{prop}[thm]{Proposition}
\newtheorem{cor}[thm]{Corollary}
\newtheorem{claim}[thm]{Claim}

\newcommand{\capsule}[1]{{\upshape\textbf{{#1}.}}}
\newcommand{\Proof}{\noindent\emph{Proof: }}


\newtheorem{makeshiftdef}{Definition}[section]
\newcommand{\Definition}[1]{ 
   \begin{makeshiftdef} 
       {\upshape {#1} } 
   \end{makeshiftdef} 
}



\newcommand{\Summary}[1]{{\small \begin{quotation} \noindent \textbf{Summary.}{#1} \end{quotation} }}
\newcommand{\Header}[1]{ \noindent \textbf{ #1 }\\ }


\title[A hypersurface whose stability group is not
determined by 2-jets]{A hypersurface in $\C^2$ whose stability group
\\is not determined by 2-jets}
\author{R. Travis Kowalski}
\address{Mathematics Department, 0112, University of California, San
Diego, La Jolla CA 92093-0112}
\email{\texttt{kowalski@math.ucsd.edu}}
\begin{thanks}{2000 {\em Mathematics Subject
Classification.} 32H12, 32V20.}\end{thanks}

\maketitle


\begin{abstract}
We give an example of a hypersurface in $\C^2$ through 0 whose stability group at 0 is determined by 3-jets, but not by jets of any lesser order.  We also examine some of the properties which the stability group of this infinite type hypersurface shares with the 3-sphere in $\C^2$.
\end{abstract}


\section{Statement of Result}

Suppose $M \subset \C^N$ is a real-analytic hypersurface
passing through the point $p$.  The \emph{stability group}
of $M$ at $p$, denoted $\mathrm{Aut}(M,p)$,  is the group
(under composition) of local automorphisms of the germ
$(M,p)$.  That is, it is the set of all invertible
biholomorphic mappings $H : \C^N \to \C^N$, defined in a
neighborhood of $p$, which fix the point $p$ and map $M$
into itself.  The stability group of $M$ at $p$ is said to
be \emph{determined by $\ell$-jets} if for every pair $H_1,
H_2 \in \mathrm{Aut}(M,p)$, we have $H_1 = H_2$ (as germs
of biholomorphisms at 0) whenever
\[
   \dd{^{|\alpha|} H_1}{Z^\alpha}(p) = \dd{^{|\alpha|}
H_2}{Z^\alpha}(p) \quad \forall \, \alpha \in \N^N, \, 0
\le |\alpha| \le \ell.
\]

Recall that a hypersurface $M \subset \C^N$ is said to be
\emph{minimal} at $p \in M$ if there exists no complex
hypersurface contained in $M$ passing through $p$.  If $M$
is real-analytic, then it is well known that this is
equivalent to being of \emph{finite type} at $p$ (in the
sense of Kohn \cite{Kohn:BoundaryBehavior} and Bloom and
Graham
\cite{BloomGraham:OnTypeConditions}).

In general, if $M \subset \C^N$ is a hypersurface of
infinite type at $p$, then its stability group at $p$ need
not be determined by jets of any finite order.  For
example, the ``flat hypersurface'' given by
\[
   M = \setofbar{ (Z_1, \dots, Z_N) \in \C^N }{ \Im \, Z_N
= 0 }
\] is of infinite type at the origin.  Moreover, any
invertible holomorphic mapping of the form
\[
   H(Z) = \big( F_1(Z), \dots F_{N-1}(Z), Z_N \big)
\] is a local automorphism of $M$.  This shows that its
stability group at 0 is not determined by $\ell$-jets for
any choice of $\ell \ge 1$.

In some sense, however, this is the most trivial example,
and for $\C^2$ in particular, it is (to the author's
knowledge) the \emph{only} such example known.  On the
other hand, there exists a large body of work concerning
the jet-determinacy of stability groups of hypersurfaces in
$\C^2$ at points of finite type.  Poincare
\cite{Poincare:LesFonctionsAnalytiques} proved that the
stability group at any point of the 3-sphere $S^3 \subset
\C^2$ is determined by 2-jets.  This was extended by Chern
and Moser \cite{ChernMoser:RealHypersurfaces}, who proved
that the stability group of a Levi-nondegenerate
hypersurface in $\C^N$ is also determined by 2-jets.  (For
more information, see the survey articles
\cite{BaouendiEbenfeltRothschild:LocalGeometricProperties}
and \cite{Vitushkin:SCVI_IV}.) More recently, Ebenfelt,
Lamel, and Zaitsev
\cite{EbenfeltLamelZaitsev:FiniteJetDetermination} have
shown that the stability group of any hypersurface of
finite type in
$\C^2$ is determined by 2-jets.

The purpose of this paper is to present an example which
shows that this result cannot be extended to nonflat
hypersurfaces in $\C^2$ of infinite type at 0 by presenting
a nonflat hypersurface $M \subset \C^2$ of infinite type
whose local automorphisms at the origin are determined by
their 3-jets, but \emph{not} by their 2-jets.

To state this result more precisely, we make one last
definition.  Let $M \subset \C^2$ be a hypersurface passing
through the origin.  A \emph{formal automorphism} of $M$ at
0 is a $\C^2$-valued invertible formal power series $H$ in
two indeterminates which vanishes at 0 and formally maps $M$
into itself.  That is, for any real-analytic local defining
function $\rho(Z,\bar{Z})$ for $M$, there exists a formal
power series $a$ in 4 indeterminates such that the
following power series identity holds:
\[
   \rho \big( H(Z), \bar{H(Z)} \big) \equiv a(Z,\bar{Z}) \,
\rho(Z,\bar{Z}).
\] The set of all such formal power series (which forms a
group under power series composition) is called the
\emph{formal stability group} of $M$ at 0, and is denoted
$\hat{\mathrm{Aut}}(M,0)$.  It is easy to see that if a
formal automorphism of $M$ converges, then it is a local
automorphism of $M$ at 0 as described above, whence it
follows that $\mathrm{Aut}(M,0) \subset
\hat{\mathrm{Aut}}(M,0)$.  We now state our main result.

\begin{thm} \label{infinite-type sphere} For the
hypersurface
\begin{equation} \label{eqn:M defined}
   M := \setofbar{ (z,w) \in \C^2 }{ |z|<1, \, \Im \, w = (
\Re
\, w )
\frac{1 - \sqrt{1 - |z|^4}}{|z|^2} },
\end{equation} 
every formal
automorphism of the germ $(M,0)$ converges.  For
$\alpha \in \C$ and $s \in \R$, let $\theta_{\alpha,s}$
be the holomorphic function defined in a neighborhood of
$0 \in \C^2$ by 
\[
   \theta_{\alpha, s}(z,w) :=  \big( 1 - 2 i \,
\bar{\alpha} \, z w - (s + i \, |\alpha|^2) w^2 \big)^{1/2},
\] where $\C \ni \zeta \mapsto (\zeta)^{1/2} \in \C$ is
 the principal branch of the square root function.    Then
the formal stability group of $M$ at 0 is given explicitly
by the following:
\begin{align}
   &\hat{\mathrm{Aut}}(M,0) = \mathrm{Aut}(M,0) 
\label{eqn:stability group}\\
   &\qquad  = \setofbar{ H^{\varepsilon,r}_{\alpha,s}(z,w)
:= \bigg( \frac{ \varepsilon(z + \alpha \, w) }{
\theta_{\alpha,s}(z,w) }, \frac{r \,
w}{\theta_{\alpha,s}(z,w) } \bigg) }{ 
       \begin{array}{c} 
          \varepsilon \in \C, \, |\varepsilon| = 1 \\
          r \in \R \setminus \{ 0 \} \\
          \alpha \in \C \\
          s \in \R
      \end{array} }. \nonumber
\end{align}
\end{thm} The proof will be given in the next section.   We
conclude this section with some remarks.\\

\noindent \textbf{Remark 1.} To the author's knowledge,
this is the first example of a nonflat hypersurface in
$\C^2$ whose stability group (at a point) is not determined
by 2-jets, or of any hypersurface in $\C^2$ whose stability
group is determined by jets of finite order, but \emph{not}
by 2-jets.  In fact, it follows from the explicit formula
above that if
\[
  \dd{^{j+k} H^{\varepsilon,r}_{\alpha,s}}{z^j \del
w^k}(0,0) = \dd{^{j+k}
H^{\varepsilon',r'}_{\alpha',s'}}{z^j \del w^k}(0,0) \quad
\forall \, j+k \le 2,
\] 
then $\varepsilon = \varepsilon'$,
$r=r'$, and $\alpha = \alpha'$, but $s$ and $s'$ are
arbitrary.  Indeed, the mappings
\begin{equation} \label{eqn:identity wannabes}
   H^{1,1}_{0,s}(z,w) = \bigg( \frac{z}{(1 - s \,
w^2)^{1/2}}, \frac{w}{(1 - s \, w^2)^{1/2}} \bigg), \quad
\sigma \in \R,
\end{equation} form a 1-parameter family of local
automorphisms of $(M,0)$ which agree with the identity
mapping up to order two, but are distinct for each
different value of $s$.\\

\noindent \textbf{Remark 2.} Observe that the hypersurface
$M$ given by equation (\ref{eqn:M defined}) is of infinite
type at 0, since it contains the nontrivial complex
hyperplane $\Sigma = \{ w = 0 \}$.  Hence, it follows that
the result of
\cite{EbenfeltLamelZaitsev:FiniteJetDetermination} does not
hold for infinite type (but nonflat) hypersurfaces in
$\C^2$.  However, it has been shown that for a particular
class of infinite type hypersurfaces (the so-called
\emph{1-infinite type} hypersurfaces, of which $M$ is an
example), stability groups are determined by jets of
\emph{some} predetermined finite order; see
\cite{EbenfeltLamelZaitsev:FiniteJetDetermination} and
\cite{Kowalski:RationalJetDependence}).\\

\noindent \textbf{Remark 3.} Since the hypersurface $M$
above is of infinite type, it is \emph{not}
biholomorphically equivalent to the 3-sphere $S^3$ in
$\C^2$.  However, the stability groups of the two
hypersurfaces have several traits in common; we point out a
few of these.
\begin{itemize}

\item It is well known that the 3-sphere in $\C^2$ is locally
biholomorphically equivalent to the hypersurface
$\setofbar{ (z,w) }{ \Im \, w = |z|^2 }$, and in these
coordinates, every (formal) local automorphism at 0 is
given by
\[
   H(z,w) = \bigg( \frac{r \, \epsilon (z + \alpha \, w)}{1
- 2i \, \bar{\alpha} z - (s + i |\alpha|^2)w}, 
                   \frac{ r^2 \, w }{1 - 2i \, \bar{\alpha}
z - (s + i |\alpha|^2)w} \bigg),
\] with $r > 0$, $\epsilon \in \C$ with $|\epsilon|=1$,
$\alpha
\in
\C$, and $s \in \R$. 
This is similar to formula given by equation
(\ref{eqn:stability group}).

\item Like the 3-sphere, the (formal) stability group of
$M$ is determined by five real parameters.

\item Like the 3-sphere, the elements of the stability
group of $M$ do not extend to a common neighborhood of 0
in $\C^2$.  That is, there exist automorphisms of the germ
$(M,0)$ whose radii of convergence are arbitrarily small. 
For example, the map $H^{1,1}_{0,s}$ given as in equation
(\ref{eqn:identity wannabes}) with $s \neq 0$ is converges
if and only if $|w| < 1/\sqrt{|s|}$, which can be made
arbitrarily small by taking $|s|$ arbitrarily large.   In
contrast, for Levi-nondegenerate hypersurfaces of $\C^2$
\emph{other than the sphere}, all local automorphisms at a
fixed point extend to a common neighborhood.

\item The stability group of $(M,0)$ forms a Lie group,
which may be identified with space $(\C - \{ 0 \}) \times
\C \times \R$ under the multiplication
\[
   ( \zeta, \alpha, s ) \cdot (\zeta' , \alpha', s') :=
   \big( \zeta \, \zeta', \alpha + \zeta \, \alpha', s + s'
- 2 \, \Im( \alpha \, \bar{\zeta} \, \bar{\alpha'} ) \big).
\] In particular, like the 3-sphere, it is noncompact,
five-dimensional, and contains a Heisenberg subgroup
(namely the subgroup defined by taking $\zeta = \zeta' =
1$).  In contrast, the stability groups of
Levi-nondegenerate hypersurfaces in $\C^2$ \emph{other than
the sphere} are compact Lie groups of dimension at most
four.

\end{itemize}

\section{Proof of Theorem \ref{infinite-type sphere}.}

We shall denote by $S^1 \subset \C$ the set of unimodular
complex numbers.  Observe that
\[
   H^{\varepsilon,r}_{\alpha, s} = \big(
H^{\varepsilon,r}_{0,0} \big) \circ \big(
H^{1,1}_{\alpha,r} \big) \quad \forall \,
(\varepsilon,r,\alpha,s) \in S^1 \times (\R \setminus
\{0\}) \times \C \times \R,
\] so to prove that $H^{\varepsilon,r}_{\alpha,s}$ is an
automorphism of $(M,0)$, it suffices to show that the
mappings $H^{\varepsilon,r}_{0,0}$ and $H^{1,1}_{\alpha,s}$
are local automorphisms of $(M,0)$.  It is obvious that the
mappings $H^{\varepsilon,r}_{0,0}$ are \emph{global}
automorphisms of $M$ for each unimodular complex number
$\varepsilon$ and nonzero real number $r$; we leave it to
the diligent reader to show that $H^{1,1}_{\alpha,r}$ is an
local automorphism of $(M,0)$ for each complex number
$\alpha$ and real number $s$.

Thus, it follows that $H^{\varepsilon,r}_{\alpha,s} \in
\mathrm{Aut}(M,0)$ for every $(\varepsilon,r,\alpha,s)\in
S^1 \times (\R \setminus \{0\}) \times \C \times \R$.  To
complete the proof, we must prove that if $H \in
\hat{\mathrm{Aut}}(M,0)$ is a formal automorphism, then $H
= H^{\varepsilon,r}_{\alpha,s}$ for some choice of
parameters $(\varepsilon,r,\alpha,s)$.   To prove this, we
introduce some new notation.  Writing $\Im \, w = (w -
\bar{w})/(2 i)$ and $\Re \, w = (w + \bar{w})/2$ in the
local defining equation (\ref{eqn:M defined}) for $M$ and
solving for $w$ yields the identity
\[
   M = \setofbar{ (z,w) }{ w = \bar{w} \, S(|z|^2) },
\] where $S$ is the real-analytic, complex-valued function
defined by
\[
   \R \supset (-1,1) \in t \mapsto S(t) := i \, t + \sqrt{1
- t^2} \in \C.
\] Recall that $H \in \hat{\mathrm{Aut}}(M,0)$ means that
$H = (H^1,H^2)$ is a $\C^2$-valued formal power series
which vanishes at 0, has nonvanishing Jacobian at 0, and
satisfies the identity
\begin{equation} \label{eqn:self map}
   H_2 \big( z, \tau \, S(z \chi) \big) \equiv
   \bar{H_2}(\chi,\tau) \, S \big( H_1 \big(z, \tau \, S(z
\chi) \big) \, \bar{H_1}(\chi,\tau) \big),
\end{equation} where $\bar{H_j}$ denotes the power series
obtained by replacing the Taylor coefficients of $H_j$ by
their complex conjugates.  Observe that if we set $\chi =
\tau = 0$ in (\ref{eqn:self map}), we obtain
\[
    H_2(z,0) = \bar{H_2}(0,0) \, S \big( H_2(z,0) \,
\bar{H_2}(0,0) \big) = 0,
\] since $H(0,0) = 0$.  Hence, we can write
\begin{equation} \label{eqn:H in f,g}
   H(z,w) = \big( f(z,w), w \, g(z,w) \big)
\end{equation} with $f(0,0) = 0$ and $f_z(0,0) \cdot g(0,0)
\neq 0$.  Substituting this into (\ref{eqn:self map}) and
cancelling a common $\tau$ from both sides yields the
identity
\begin{equation} \label{eqn:self map 2}
   S(z \chi) g \big( z, \tau \, S(z \chi) \big) \equiv
   \bar{g}(\chi, \tau) \, S \big( f \big( z, \tau \, S(z
\chi) \big) \, \bar{f}(\chi,\tau) \big).
\end{equation} Finally, for convenience, we shall formally
expand the power series $f$ and $g$ as
\begin{equation} \label{eqn:expanded f and g}
   f(z,w) = \sum_{n=0}^\infty \frac{f_n(z)}{n!} w^n, \quad
   g(z,w) = \sum_{n=0}^\infty \frac{g_n(z)}{n!} w^n,
\end{equation} and shall write
\begin{equation} \label{eqn:barred derivatives}
   a_n^j := \bar{ f_n^{(j)}(0) }, \quad b_n^j := \bar{
g_n^{(j)}(0) }, \quad n, j \ge 0.
\end{equation} We now state the main lemma which will
complete the proof of Theorem \ref{infinite-type sphere}.

\begin{lem} \label{finite det} Let $M$ be the hypersurface
defined in Theorem \ref{infinite-type sphere}.  Suppose
that $H \in \hat{\mathrm{Aut}}(M,0)$, and write $H$ as in
equation (\ref{eqn:H in f,g}).  Then for every $n \ge 0$,
there exists a $\C^2$-valued polynomial $R_n$ in eight
indeterminates, depending only on $M$ and not the formal
map $H$, such that
\[
   \big( f_n(z), g_n(z) \big) = R_n \bigg( z,
\frac{1}{a_0^1}, \frac{1}{b_0^0}, a_0^1, b_0^0, a_1^0,
\bar{a_1^0}, \Re \, b_2^0 \bigg).
\] Moreover, we have $ a_0^1 \in S^1$ and $b_0^0 \in \R
\setminus \{ 0 \}$.
\end{lem}

To see that Lemma \ref{finite det} completes the proof of
Theorem \ref{infinite-type sphere}, fix a formal
automorphism $H \in \hat{\mathrm{Aut}}(M,0)$.  Lemma
\ref{finite det} implies $H$ is \emph{uniquely determined}
by its values $a_0^1$, $b_0^0$, $a_1^0$ and $\Re \,
b_2^0$.  Define
\[
   \varepsilon := \bar{a_0^1} \in S^1, \quad
   r := \bar{b_0^0} \in \R \setminus \{ 0 \}, \quad
   \alpha := \frac{ (\bar{a_1^0}) }{ (\bar{a_0^1}) } \in
\C, \quad
   s := \frac{ \Re( \bar{b_2^0} ) }{ (\bar{b_0^0}) } \in \R,
\] and define $\tilde{H} := H^{\varepsilon,r}_{\alpha,s}$. 
Define the corresponding (barred) derivatives
$\tilde{a}_n^j$ and $\tilde{b}_n^j$ for $\tilde{H}$ as in
equations (\ref{eqn:H in f,g}), (\ref{eqn:expanded f and
g}), and (\ref{eqn:barred derivatives}).  It follows from a
simple calculation that $a_0^1 = \tilde{a}_0^1$, $b_0^0 =
\tilde{b}_0^0$, $a_1^0 = \tilde{a}_0^1$, and $\Re \, b_2^0
= \Re \, \tilde{b}_2^0$, whence $H = \tilde{H}$ by
uniqueness, and the proof of the Theorem is complete. 
Hence, we need only prove the lemma.\\

\noindent \emph{Proof of Lemma \ref{finite det}:} We
proceed by induction.  For convenience, we shall set
\[
   \lambda_0 := \bigg( \frac{1}{a_0^1}, \frac{1}{b_0^0},
a_0^1, b_0^0, a_1^0, \bar{a_1^0}, \Re \, b_2^0 \bigg) \in
\C^7.
\] For any formal power series $H$ of the form (\ref{eqn:H
in f,g}), define
\[
   \Phi^H(z,\chi,\tau) := - S(z \chi) g \big( z, \tau \,
S(z \chi) \big) +
   \bar{g}(\chi, \tau) \, S \big( f \big( z, \tau \, S(z
\chi) \big) \, \bar{f}(\chi,\tau) \big).
\] By (\ref{eqn:self map 2}), it follows that an invertible
power series $H$ is a formal automorphism of $(M,0)$ if and
only if $\Phi^H \equiv 0$.  The basic algorithm of the
proof is as follows: given a formal automorphism $H$, at
the $n$-th step of the induction, we
\begin{itemize}
   \item Calculate $\Phi^H_{\tau^n}(z, \chi, 0)$.
   \item Solve $\Phi^H_{\tau^n}(z,0,0) = 0$ to obtain an
explicit formula for $g_n(z)$ as a polynomial (independent
of the mapping $H$) in $(z,a_n^0, a_n^1, b_n^0, b_n^1,
\lambda_0) \in \C^{12}$.
   \item Solve $\Phi^H_{\chi \, \tau^n}(z,0,0) = 0$ to
obtain an explicit formula for $f_n(z)$, similarly
expressed.
   \item Substitute these formulas (and their complex
conjugates) into the identity $\Phi^H_{\tau^n}(z,\chi,0) =
0$ and differentiate this repeatedly in $z$ and $\chi$  to
express $(a_n^0, a_n^1, b_n^0, b_n^1)$ as a polynomial in
$\lambda_0$.
\end{itemize} In the algorithm above, we have used the
usual subscript notation to denote partial derivatives,
i.~e.
\[
   \Phi^H_{z^j \chi^k \tau^\ell}(z,\chi,\tau) :=
\dd{^{j+k+\ell} \Phi^H}{z^j \del \chi^k \del
\tau^\ell}(z,\chi,\tau).
\] We now fill in the details.  Fix an automorphism $H$.

\emph{The case $n=0$.}  Setting $\Phi^H(z,0,0) = 0$, we
obtain $g_0(z) = \bar{g_0}(0) = b_0^0$, from which it
follows that $b_0^0$ is real and, since $H$ is invertible,
nonzero.  Thus, we have
\begin{equation} \label{eqn:g_0}
   g_0(z) = \bar{g_0}(\chi) = b_0^0 \in \R \setminus \{ 0
\}.
\end{equation} Setting $\Phi^H_\chi(z,0,0) = 0$ and using
(\ref{eqn:g_0}), we find $f_0(z) = z/a_0^1$.  From this, it
follows that $\bar{a_0^1} = 1/a_0^1$, so $a_0^1$ is
necessarily unimodular.  Thus, we have
\begin{equation} \label{eqn:f_0}
   f_0(z) = \frac{z}{a_0^1}, \quad \bar{f_0}(\chi) = a_0^1
\, \chi, \quad
   a_0^1 \in S^1,
\end{equation} which completes the base step of the
induction.

\emph{The case $n=1$.} Using the identity
$\Phi^H_\tau(z,\chi,0) = 0$ as indicated above and
substituting in the formulas (\ref{eqn:g_0}) and
(\ref{eqn:f_0}) as needed, we find
\[
   f_1(z) =  \frac{i \, a_1^0}{(a_0^1)^2} \, z^2
     + \bigg( \frac{ b_1^0}{a_0^1 \, b_0^0} -
\frac{a_1^1}{(a_0^1)^2} \bigg)z 
     + \frac{i \, b_1^1}{a_0^1 \, b_0^0}, \quad
   g_1(z) =  \frac{i \, b_0^0 \, a_1^0}{a_0^1} \, z + b_1^0.
\] Conjugating these, we obtain
\[
   \bar{f_1}(\chi)  = \frac{a_0^1 \, b_1^1}{b_0^0} \, \chi^2
     + a_1^1 \, \chi + a_1^0 ,\quad
   \bar{g_1}(\chi)  = b_1^1 \, \chi + b_1^0.
\] Using these formulas, it follows that
\begin{equation} \label{eqn:n=1 reduction}
   \left( \begin{array}{c} 0 \\ 0 \end{array} \right)
   = \left( \begin{array}{c} 
       \Phi^H_{z^2 \chi^2 \tau}(0,0,0) \\
       \Phi^H_{z^3 \chi^3 \tau}(0,0,0)
     \end{array} \right)
  = \left( \begin{array}{c c}
         4 \frac{ b_0^0 }{ a_0^1 } & - 2 \\
      54 i \frac{ b_0^0 }{ a_0^1 } &  - 18 i
    \end{array} \right)
    \left( \begin{array}{c}
      a_1^1 \\ b_1^0
    \end{array} \right).
\end{equation} Since the $2\times2$ matrix on the righthand
side of equation (\ref{eqn:n=1 reduction}) is invertible,
it follows from equation (\ref{eqn:n=1 reduction}) that
$a_1^1 = b_1^0 = 0$.  Moreover, equating $\bar{a_1^0} =
f_1(0)$ yields $b_1^1 = - i \, a_0^1 \, b_0^0 \,
\bar{a_1^0}$.  Hence, we have
\[
   f_1(z) =  \frac{i \, a_1^0}{(a_0^1)^2} \, z^2 +
\bar{a_1^0}, \quad
   g_1(z) =  \frac{i \, b_0^0 \, a_1^0}{a_0^1} \, z,
\] which completes the induction at this step.

\emph{The case $n=2$.} Using the identity
$\Phi^H_{\tau^2}(z,\chi,0) = 0$ as above, we find
\begin{align*}
   f_2(z) &= -   \frac{3(a_1^0)^2}{(a_0^1)^3} \, z^3
    +   \frac{2i \, a_2^0}{(a_0^1)^2} \, z^2
    + \bigg( \frac{2 i \, a_1^0 \, \bar{a_1^0} }{a_0^1} 
       + \frac{2 b_2^0}{a_0^1 \, b_0^0} -
\frac{a_2^1}{(a_0^1)^2} \bigg) z
    + \frac{i \, b_2^1}{a_0^1 \, b_0^0}, \\
   g_2(z) &= -  \frac{3 b_0^0 \, (a_1^0)^2}{(a_0^1)^2} \,
z^2 
    +  \frac{i \, b_0^0 \, a_2^0}{a_0^1} \, z
    +  b_2^0 + 2 i \, b_0^0 \, a_1^0 \, \bar{a_1^0}.
\end{align*} Conjugating as above and substituting into
$\Phi^H_{\tau^2}(z,\chi,0) = 0$, the relations 
\[
   \Phi^H_{z^2 \chi^2 \tau^2}(0,0,0)=0, \quad 
   \Phi^H_{z^3 \chi^2 \tau^2}(0,0,0)=0, \quad
   \Phi^H_{z^2 \chi^3 \tau^2}(0,0,0)=0 
\]  yield $a_2^0 = b_2^1 = 0$ and $a_2^1 = (a_0^1 \,
b_2^0)/b_0^0 - 2 i \, a_0^1 \, a_1^0 \, \bar{a_1^0}$. 
Similarly, equating $\bar{b_2^0} = g_2(0)$, we find
\[
   \Im \, b_2^0 = \frac{ b_2^0 - \bar{b_2^0} }{2 i} = -
\frac{i \, a_1^0 \,b_1^1}{a_0^1} = - b_0^0 \, a_1^0
\, \bar{a_1^0}.
\] Under these substitutions, we have
\begin{align*}
   f_2(z) &= - \frac{3(a_1^0)^2}{(a_0^1)^2} \, z^3
     + \bigg( \frac{\Re \, b_2^0}{a_0^1 \, b_0^0} 
         + \frac{3 i \, a_1^0 \, \bar{a_1^0}}{a_0^1} \bigg)
z, \\
   g_2(z) &= -  \frac{3 b_0^0 (a_1^0)^2}{(a_0^1)^2} \, z^2
     + \Re \, b_2^0 + i \, b_0^0 \, a_1^0 \, \bar{a_1^0} ,
\end{align*} which completes the induction at this step.

\emph{The general inductive step.}  Assume now that the
lemma holds up to some $n - 1 \ge 2$; we prove it for $n$. 
The Chain Rule implies
\begin{align}
   &\Phi^H_{\tau^n}(z,\chi,0)
   = - S(z \chi)^{n+1} g_n(z) 
      + S \big( z \chi \big) \bar{g_n}(\chi)  \nonumber \\
   &\qquad  + b_0^0 \, S'\big( z \chi \big) 
         \bigg( a_0^1 \, \chi \, S(z \chi)^n f_n(z) +
\frac{z}{a_0^1} \, \bar{f_n}(\chi) \bigg)
\label{eqn:general step}\\
   &\qquad  + P^n \bigg( \big( S^{(j)}(z \chi)
\big)_{j=0}^n, \big( f_j(z), g_j(z), \bar{f_j}(\chi),
\bar{g_j}(\chi) \big)_{j=0}^{n-1} \bigg), \nonumber
\end{align} where $P_n$ is a complex-valued polynomial (in
$5n+1$ indeterminates) which is independent of the mapping
$H$.  By the inductive hypothesis (and its conjugation), we
may rewrite the last term in equation (\ref{eqn:general
step}) as
\[
   Q^n \bigg (z, \chi, \lambda_0, \big( S^{(j)}(z \chi)
\big)_{j=0}^n \bigg),
\] where $Q^n$ is complex polynomial in $n+10$
indeterminates, independent of $H$.  Proceeding as above,
we find
\begin{align} 
   f_n(z) &=  \frac{i \, n \, a_n^0}{(a_0^1)^2} \, z^2
     + \bigg( \frac{n \, b_n^0}{a_0^1 \, b_0^0} 
        - \frac{a_n^1}{(a_0^1)^2} \bigg)z 
     + \frac{i \, b_n^1}{a_0^1 \, b_0^0} + p_n (z,
\lambda_0),  \label{eqn:f_n} \\
   g_n(z) &=  \frac{i \, b_0^0 \, a_n^0}{a_0^1} z + b_n^0 
     + q^n(z,\lambda_0), \label{eqn:g_n}
\end{align} where $p_n, \, q_n$ are complex polynomials in
8 indeterminates, independent of $H$.  Substituting these
and their conjugates into the identity
$\Phi^H_{\tau^n}(z,\chi,0) = 0$ and then computing
$\Phi^H_{z^j \chi^k \tau^n}(0,0,0)$ for $j,k = 2,3$ yields
a $4 \times 4$ system of equations of the form
\begin{equation} \label{eqn:general reduction}
   A_n \cdot \big( a_n^0,a_n^1,b_n^0,b_n^1 )^t =
B_n(\lambda_0),
\end{equation} where $B_n$ is a $\C^4$-valued polynomial in
$\lambda_0$, and $A_n$ is the $4 \times 4$ matrix given by
\[
   \left( \begin{array}{cccc}
     0 & 4n \frac{b_0^0}{a_0^1} & -2 n^2 & 0 \\
     -6 i (n^2-1) \frac{b_0^0}{a_0^1} & 0 & 0 & 0 \\
     0 & 0 & 0 & 6 (n^2-1) \\
     0 & 18 i (n^2+2n) \frac{b_0^0}{a_0^1} & -
6i(2n^3+3n^2-2n) & 0
   \end{array} \right).
\] Observe that
\[
   \det(A_n) = - \frac{432(b_0^0)^2}{(a_0^1)^2}
    (n-2)(n-1)^2 n^2 (n+1)^2 (n+2),
\] which is nonzero for $n \ge 3$, whence $A_n$ is
invertible.  By Craemer's Rule, it follows that $A_n^{-1}$
is a $4 \times 4$ matrix whose entries are polynomial in
$(a_0^1,b_0^0)$ and their reciprocals (and so in particular
are polynomial in $\lambda_0$).  Thus, equation
(\ref{eqn:general reduction}) implies that $\big(
a_n^0,a_n^1,b_n^0,b_n^1 )$ is a polynomial in $\lambda_0$. 
Substituting this into equations (\ref{eqn:f_n}) and
(\ref{eqn:g_n}) completes the induction. \qed

\bibliographystyle{alpha}

\begin{thebibliography}{BER00}

\bibitem[BER00]{BaouendiEbenfeltRothschild:LocalGeometricProperties}
M.~S. Baouendi, P.~Ebenfelt, and L.~P. Rothschild.
\newblock Local geometric properties of real submanifolds in complex space.
\newblock {\em Bull. Amer. Math. Soc. (N.S.)}, 37(3):309--336 (electronic),
  2000.

\bibitem[BG77]{BloomGraham:OnTypeConditions}
T.~Bloom and I.~Graham.
\newblock On ``type'' conditions for generic real submanifolds of
  $\mathbb{C}^n$.
\newblock {\em Invent. Math.}, 40(3):217--243, 1977.

\bibitem[CM74]{ChernMoser:RealHypersurfaces}
S.~S. Chern and J.~K. Moser.
\newblock Real hypersurfaces in complex manifolds.
\newblock {\em Acta Math.}, 133:219--271, 1974.

\bibitem[ELZ00]{EbenfeltLamelZaitsev:FiniteJetDetermination}
P.~Ebenfelt, B.~Lamel, and D.~Zaitsev.
\newblock Finite jet determination of local analytic {CR} automorphisms and
  their parametrization by 2-jets in the finite type case. Preprint
2001.
\newblock \texttt{http://arXiv.org/abs/math.CV/0107013}

\bibitem[Koh72]{Kohn:BoundaryBehavior}
J.~J. Kohn.
\newblock Boundary behavior of $\bar{\del}$ on weakly pseudo-convex manifolds
  of dimension two.
\newblock {\em J. Differential Geometry}, 6:523--542, 1972.

\bibitem[Kow01]{Kowalski:RationalJetDependence}
R.~T. Kowalski.
\newblock Rational jet dependence of formal equivalences between real-analytic
  hypersurfaces in $\mathbb{C}^2$.
\newblock Preprint, 2001.

\bibitem[Poi07]{Poincare:LesFonctionsAnalytiques}
H.~Poincar\'e.
\newblock Les fonctions analytiques de deux variables et la repr\'esentation
  conforme.
\newblock {\em Rend. Circ. Mat. Palermo, II. Ser.}, 23:544--547, 1907.

\bibitem[Vit90]{Vitushkin:SCVI_IV}
A.~G. Vitushkin.
\newblock {\em Several complex variables. {I}}.
\newblock Springer-Verlag, Berlin, 1990.
\newblock Translation by P. M. Gauthier.

\end{thebibliography}

\end{document}